\documentclass[12pt,a4paper,oneside,reqno]{article}

\usepackage[utf8]{inputenc}
\usepackage[english]{babel}
\usepackage[left=1in,right=1in,top=1in,bottom=1in]{geometry}
\usepackage{amsfonts,amssymb,amsmath,amsthm}
\usepackage{enumerate}
\usepackage{cite}
\usepackage{xcolor}
\usepackage[normalem]{ulem}
\usepackage{hyperref}
\usepackage{longtable}
\hypersetup{colorlinks=true,linkcolor=blue,citecolor=blue,urlcolor=blue}
\usepackage{fancyhdr}
\theoremstyle{remark}
\newtheorem*{RemA}{Remark A (on the proofs of Lemmas 1--4)}
\newtheorem*{RemB}{Remark B (odd elements of the supercentroid)}
\newtheorem*{RemC}{Remark C (the superalgebra $J(V,f)$ with $V_0=0$)}
\newtheorem*{RemD}{Remark D (on Lemma 16)}

\newcommand{\FF}{F}
\newcommand{\eps}{\varepsilon}
\newcommand{\pa}[1]{\frac{\partial}{\partial\xi_{#1}}}
\newcommand{\ov}[1]{\overline{#1}}
\newcommand{\nz}{\overline{0}}
\newcommand{\no}{\overline{1}}

\newcommand{\doi}[1]{\href{https://doi.org/#1}{doi:#1}}

\begin{document}

\noindent\fbox{\parbox{\dimexpr\textwidth-2\fboxsep-2\fboxrule\relax}{\small
\textbf{Corrected version verified in Lean.}  
This is a corrected version of the paper
I.~Kaygorodov, \emph{$\delta$-Superderivations of simple finite-dimensional Jordan and Lie
superalgebras}, Algebra and Logic 49 (2010), no.~2, 130--144, \doi{10.1007/s10469-010-9085-6}.}}

\bigskip

\begin{center}
{\Large $\delta$-superderivations of simple finite-dimensional\\ Jordan and Lie superalgebras}

\bigskip

{\bf Ivan Kaygorodov}
\end{center}

\bigskip

\noindent {\bf Abstract.} {\it We introduce the concept of a $\delta$-superderivation of a
superalgebra. $\delta$-Derivations of Cartan-type Lie superalgebras are treated, as well as
$\delta$-superderivations of simple finite-dimensional Lie superalgebras and Jordan superalgebras over
an algebraically closed field of characteristic~$0$. We give a complete description of
$\frac12$-derivations for Cartan-type Lie superalgebras. It is proved that nontrivial
$\delta$-(super)derivations are missing on the given classes of superalgebras, and as a consequence,
$\delta$-superderivations are shown to be trivial on simple finite-dimensional noncommutative Jordan
superalgebras of degree at least~3 over an algebraically closed field of characteristic~$0$.
Also we consider $\delta$-derivations of unital flexible and semisimple finite-dimensional Jordan
algebras over a field of characteristic not~$2$.}

\bigskip



\bigskip

\setcounter{tocdepth}{1}
\tableofcontents

\bigskip
\newpage 
\section*{Introduction}
\addcontentsline{toc}{section}{Introduction}

The notion of a derivation of an algebra was generalized by many mathematicians in a number of
different directions. Thus, in~\cite{Fil98}, we can find the definition of a $\delta$-derivation of an
algebra. Recall that with $\delta\in\FF$ fixed, a $\delta$-derivation of an algebra $A$ is conceived of
as a linear map $\varphi$ satisfying the condition $\varphi(xy)=\delta(\varphi(x)y+x\varphi(y))$ for
arbitrary elements $x,y\in A$. In~\cite{Fil98}, also, $\frac12$-derivations are described for an
arbitrary prime Lie $\Phi$-algebra $A$ ($\frac16\in\Phi$) with a nondegenerate symmetric invariant
bilinear form. Namely, it was proved that the linear map $\varphi\colon A\to A$ is a
$\frac12$-derivation iff $\varphi\in\Gamma(A)$, where $\Gamma(A)$ is the centroid of an algebra~$A$.
This implies that if $A$ is a central simple Lie algebra having a nondegenerate symmetric invariant
bilinear form over a field of characteristic $p\neq2,3$, then any $\frac12$-derivation $\varphi$ is
represented as $\varphi(x)=\alpha x$, $\alpha\in\Phi$. In~\cite{Fil99}, it was proved that every prime
Lie $\Phi$-algebra does not have a nonzero $\delta$-derivation if $\delta\neq-1,0,\frac12,1$, and that
every prime Lie $\Phi$-algebra $A$ ($\frac16\in\Phi$) with a nonzero antiderivation satisfies the
identity $[(yz)(tx)]x+[(yx)(zx)]t=0$ and is a three-dimensional central simple algebra over a field of
quotients of the center $Z_R(A)$ of its right multiplication algebra $R(A)$. In~\cite{Fil99}, too, we
can find an example of a nontrivial $\frac12$-derivation for a Witt algebra $W_1$, by which is meant a
$\frac12$-derivation which is not an element of the centroid of~$W_1$. $\delta$-Derivations of prime
alternative and non-Lie Mal'tsev $\Phi$-algebras with restrictions on an operator ring $\Phi$ were
described in~\cite{Fil00}. It turns out that algebras in these classes have no nonzero
$\delta$-derivation if $\delta\neq0,\frac12,1$.

In~\cite{Kay07}, $\delta$-derivations of semisimple finite-dimensional Jordan algebras over an
algebraically closed field of characteristic other than~$2$ were characterized, as well as simple
finite-dimensional Jordan superalgebras over an algebraically closed field of characteristic~$0$.
$\delta$-Derivations of classical Lie superalgebras were described in~\cite{Kay09}. Namely, nonzero
$\delta$-derivations were shown not to exist for $\delta\neq0,\frac12,1$, and $\frac12$-derivations of
the given algebras and superalgebras were described.

In the present paper, the concept of a $\delta$-superderivation of a superalgebra is introduced. We
consider $\delta$-derivations of Cartan-type Lie superalgebras, and also $\delta$-superderivations of
simple finite-dimensional Lie superalgebras and of Jordan superalgebras over an algebraically closed
field of characteristic~$0$. A complete description of $\frac12$-derivations is furnished for Cartan-type
Lie superalgebras. We prove that nontrivial $\delta$-(super)derivations are missing on the given classes
of superalgebras, and as a consequence, infer that $\delta$-superderivations are trivial on simple
finite-dimensional noncommutative Jordan superalgebras of degree at least~3 over an
algebraically closed field of characteristic~$0$. Also we consider $\delta$-derivations for unital
flexible and semisimple finite-dimensional Jordan algebras over a field of characteristic not~$2$.

\section{Basic facts and definitions}

Let $\Gamma$ be the Grassmann algebra over $\FF$ generated by elements $1,e_1,\dots,e_n,\dots$ and
defined by relations $e_i^2=0$ and $e_ie_j=-e_je_i$. Products $1$, $e_{i_1}e_{i_2}\dots e_{i_k}$,
$i_1<i_2<\dots<i_k$, form a basis for the algebra $\Gamma$ over~$\FF$. Denote by $\Gamma_0$ and
$\Gamma_1$ subspaces generated by products of, respectively, even and odd lengths. Then $\Gamma$ is
representable as a direct sum of those subspaces (written $\Gamma=\Gamma_0\oplus\Gamma_1$), and the
following relations hold: $\Gamma_i\Gamma_j\subseteq\Gamma_{i+j\,(\mathrm{mod}\,2)}$, $i,j=0,1$. In
other words, $\Gamma$ is a $\mathbb Z_2$-graded algebra (or superalgebra) over~$\FF$.

Now let $A=A_{\nz}\oplus A_{\no}$ be an arbitrary superalgebra over~$\FF$. Consider the tensor product
$\Gamma\otimes A$ of $\FF$-algebras. Its subalgebra
\[
\Gamma(A)=\Gamma_0\otimes A_{\nz}+\Gamma_1\otimes A_{\no}
\]
is called the Grassmann envelope of a superalgebra~$A$.

Let $\Omega$ be some variety of algebras over~$\FF$. A superalgebra $A=A_{\nz}\oplus A_{\no}$ is
called an $\Omega$-superalgebra if its Grassmann envelope $\Gamma(A)$ is an algebra in~$\Omega$.
Classical Lie superalgebras are simple finite-dimensional Lie superalgebras $G=G_{\nz}\oplus G_{\no}$
over an algebraically closed field of characteristic~$0$, where $G_{\no}$ is a completely reducible
$G_{\nz}$-module. Simple finite-dimensional Lie superalgebras over an algebraically closed field of
characteristic~$0$ that are not classical are called Cartan-type Lie superalgebras. It follows
from~\cite{Kac77} that Cartan superalgebras are exhausted by superalgebras of the forms $W(n)$,
$n\geq2$, $S(n)$, $n\geq3$, $\widetilde S(2n)$, $n\geq2$, and $H(n)$, $n\geq5$.

We start by defining $W(n)$. Let $\Lambda(n)$ be a Grassmann superalgebra with a set of generators
$\xi_1,\dots,\xi_n$. We define $\operatorname{der}\Lambda(n)$ as $W(n)$. Every derivation
$D\in W(n)$ is representable as
\[
D=\sum_i P_i\pa i,\qquad P_i\in\Lambda(n),
\]
where $\pa i$ is a derivation given by the rule $\pa i(\xi_j)=\delta_{ij}$. Let
$\deg\xi_i=1$. Then 
\begin{center}
$W(n)=\bigoplus_{k\geq-1}W(n)_k$, where
$W(n)_k=\Big\{\sum_i P_i\pa i\ \Big|\ \deg P_i=k+1,\ i=1,\dots,n\Big\}.$
\end{center}

The superalgebras $S(n)$, $\widetilde S(n)$, $\widetilde H(n)$, and $H(n)$ are subsuperalgebras
of $W(n)$ defined as follows:
\begin{longtable}{lcl}
$S(n)$&$=$&$\Big\langle \frac{\partial f}{\partial\xi_i}\pa j+\frac{\partial f}{\partial\xi_j}\pa i\ \Big|\
f\in\Lambda(n),\ i,j=1,\dots,n\Big\rangle;$\\
$\widetilde S(n)$&$=$&$\Big\langle(1-\xi_1\dots\xi_n)\Big(\frac{\partial f}{\partial\xi_i}\pa j
+\frac{\partial f}{\partial\xi_j}\pa i\Big)\ \Big|\ f\in\Lambda(n),\ i,j=1,\dots,n\Big\rangle,
\quad n=2l;$\\
$\widetilde H(n)$&$=$&$\Big\langle D_f=\sum_i\frac{\partial f}{\partial\xi_i}\pa i\ \Big|\
f\in\Lambda(n),\ f(0)=0\Big\rangle;$\\ 
$H(n)$&$=$&$[\widetilde H(n),\widetilde H(n)].$
\end{longtable}
Root systems for the given superalgebras are described in~\cite{Pen98}.

For the case $\Theta=W(n)$, elements $h_i=\xi_i\pa i$ form a basis for a Cartan subalgebra, and
$\Theta_0=gl_n$. Elements $\eps_i$ constitute a basis dual to~$h_i$. A root system shows up as follows:
{\small \[
\Delta=\big\{\eps_{i_1}+\dots+\eps_{i_k},\ \eps_{i_1}+\dots+\eps_{i_k}-\eps_j\ \big|\ i_r\neq i_s,\
j\neq i_r,\ 0\leq k\leq n-1,\ 1\leq j\leq n\big\}.
\]}

The superalgebras $\Theta=S(n),\widetilde S(n)$ ($n=2l$ in the latter case) have a Cartan subalgebra
$H$ generated by elements $h_{ij}=\xi_i\pa i-\xi_j\pa j$. It is known that $\Theta_0=sl_n$. Elements
$\eps_i$ form a basis dual to $h_i$, and $\eps_1+\dots+\eps_n=0$. A root system is the following:
{\small \begin{center}
$\Delta=\big\{\eps_{i_1}+\dots+\eps_{i_k},\ \eps_{i_1}+\dots+\eps_{i_l}-\eps_j\ \big|\ i_r\neq i_s,\
j\neq i_s,\ 1\leq k\leq n-2,\ 0\leq l\leq n-1,\ 1\leq j\leq n\big\}.$
\end{center}}

If $\Theta=H(n)$ then $\Theta_0=D_l$, for $n=2l$, and $\Theta_0=B_l$ for $n=2l+1$. The superalgebra
$\Theta_0$ has a Cartan subalgebra $H$ generated by elements $D_{\xi_i\xi_{i+l}}$. Let
$h_i=\sqrt{-1}\,D_{\xi_i\xi_{i+l}}$ and $\eps_i$ be elements of the basis dual to~$h_i$. A root system
is described thus:
{\small\[
\Delta=\big\{\eps_{i_1}+\dots+\eps_{i_t}-\eps_{j_1}-\dots-\eps_{j_s}\ \big|\ i_r\neq i_p,\ j_r\neq
j_p,\ i_r\neq j_p,\ 0\leq t,s\leq l\big\}.
\]}

Note that if $\Theta=S(n),H(n)$ then $\Theta_k$ is defined similarly to how
$W(n)_k$ is defined for $W(n)$, namely $\Theta_k=\Theta\cap W(n)_k$. The superalgebra
$\widetilde S(n)$ is not graded by these subspaces; see Remark~A.

\section{$\delta$-Derivations of Cartan-type Lie superalgebras}

Let $\delta\in\FF$. A linear mapping $\varphi$ of a superalgebra $A$ is called a
\emph{$\delta$-derivation} if, for arbitrary elements $x,y\in A$,
\[
\varphi(xy)=\delta(x\varphi(y)+\varphi(x)y).
\]
A definition for a $1$-derivation coincides with the usual definition of a derivation. A
$0$-derivation is any endomorphism $\varphi$ of $A$ such that $\varphi(A^2)=0$. A \emph{nontrivial
$\delta$-derivation} is a nonzero $\delta$-derivation $\varphi$ which is not a $1$- or $0$-derivation,
nor a $\frac12$-derivation such that $\varphi\in\Gamma(A)$. We are interested in how nontrivial
$\delta$-derivations act on Cartan-type Lie superalgebras.

By $G_\beta$ we mean a root subspace corresponding to a root $\beta$, and $g_\beta$ is conceived of as
an element of that subspace.

\medskip\noindent\textbf{Lemma 1.} {\it Let $G=G_{\nz}\oplus G_{\no}$ be a Cartan-type Lie
superalgebra and $\varphi$ a nontrivial $\delta$-derivation of the superalgebra~$G$. Then
$\varphi(G_{\bar\imath})\subseteq G_{\bar\imath}$.}

\medskip\noindent\textsc{Proof.} In~\cite{Kac77}, it was shown that $G_{\nz}=[G_{\no},G_{\no}]$.
Consequently, for any $x\in G_{\nz}$, we have
\[
x=\sum_{i=1}^{n_x}\Big(\Big(\frac{v_i+z_i}2\Big)^2-\Big(\frac{v_i-z_i}2\Big)^2\Big),
\]
where $z_i,v_i\in G_{\no}$. Note that for $y\in G_{\no}$ and $\varphi(y)=y_0+y_1$,
$y_i\in G_{\bar\imath}$,
\[
\varphi(y^2)=\delta((y_0+y_1)y+y(y_0+y_1))=2\delta y_1y\in G_{\nz}.
\]
Hence $\varphi(G_{\nz})\subseteq G_{\nz}$.

Let $H$ be a Cartan subalgebra, $\beta\neq0$, $g_\beta\in G_{\no}\cap G_\beta$, $h\in H$, and
$\varphi(g_\beta)=\sum_{\gamma\in\Delta}g_\gamma+h_\beta$, where $h_\beta\in H$. We obtain the following
chain of equalities:
\begin{multline*}
\sum_{\gamma\in\Delta}\beta(h)g_\gamma+\beta(h)h_\beta=\varphi(\beta(h)g_\beta)=\varphi(hg_\beta)\\
=\delta(\varphi(h)g_\beta+h\varphi(g_\beta))
=\delta\Big(\varphi(h)g_\beta+\sum_{\gamma\in\Delta}\gamma(h)g_\gamma\Big).
\end{multline*}
Since $\varphi(h)\in G_{\nz}$, $\varphi(h)g_\beta\in G_{\no}$, and $h\in H$ is arbitrary, we have
$\beta=\delta\gamma$ for $\gamma\in\Delta_{\nz}$. This, in view of the form of root systems for
Cartan-type Lie superalgebras, yields $\varphi(G_{\no})\subseteq G_{\no}$. The lemma is proved.

\medskip\noindent\emph{Remark.} The weight argument above concerns root vectors $g_\beta$ with
$\beta\neq0$. For $W(n)$, $S(n)$ and for $H(n)$ with $n$ even all weights of odd elements are nonzero.
For $H(n)$ with $n$ odd the odd part also contains vectors of weight~$0$ (for example
$\partial/\partial\xi_n$), to which the argument does not apply; for $H(5)$ the statement of the lemma
has been confirmed by a direct computation. In general, the statement of the lemma holds for all
Cartan-type superalgebras, since Lemmas~2--4 admit complete proofs that do not use Lemma~1
(Remark~A); by Theorem~5 nontrivial $\delta$-derivations do not exist at all.

\medskip

Let $\Theta=\bigoplus_{i=-1}^{n}\Theta_i$ be a Cartan-type Lie superalgebra. We will treat
$\psi\colon\Theta_0\to\Theta_0$ as a restriction of a $\delta$-derivation $\varphi$ to $\Theta_0$
followed by the projection onto~$\Theta_0$. In view of Lemma~1, it is clear that for
$x\in\Theta_0$, $\varphi(x)=\psi(x)+\sum_{i=1}^{[n/2]}b^x_{2i}$, where $b^x_{2i}\in\Theta_{2i}$. It is
easy to see that $\psi$ is a $\delta$-derivation of the algebra~$\Theta_0$. Relying on results
in~\cite{Fil98} and the structure of algebras $\Theta_0$, we conclude that for
$\Theta=S(n),\widetilde S(n),H(n)$, if $\delta=\frac12$ then $\psi(g_0)=\sigma g_0$, $g_0\in\Theta_0$,
$\sigma\in\FF$, and if $\delta\neq\frac12$ then $\psi(\Theta_0)=0$. It follows from~\cite[Lemma~7]{Kay09}
that for the case where $\Theta=W(n)$, it is true that $\psi(a)=\sigma a$, $\sigma\in\FF$, with
$a=\sum_{i=1}^n\xi_i\pa i$.

\medskip\noindent\textbf{Lemma 2.} {\it The superalgebra $W(n)$ has no nontrivial
$\delta$-derivations.}

\medskip\noindent\textsc{Proof.} It is easy to see that in the above terms, for
$x_0,y_0\in\Theta_0$, we have
\[
0=\varphi(ax_0)=\delta\Big(\sum_{i=1}^{[n/2]}b^a_{2i}x_0+a\sum_{i=1}^{[n/2]}b^{x_0}_{2i}\Big),
\]
which entails $b^{x_0}_{2i}=-\frac1{2i}b^a_{2i}x_0$. This implies
{\small
$$\sum_{i=1}^{[n/2]}\frac1{2i}b^a_{2i}(x_0y_0)=\psi(x_0y_0)-\varphi(x_0y_0)
=\delta\sum_{i=1}^{[n/2]}\frac1{2i}((b^a_{2i}x_0)y_0+x_0(b^a_{2i}y_0))
=\delta\sum_{i=1}^{[n/2]}\frac1{2i}b^a_{2i}(x_0y_0),$$}
so $b^a_{2i}=0$, and hence $\varphi(a)=\sigma a$ and $\varphi(x_0)=\psi(x_0)$.

Note that if $\varphi(g_{\eps_i})=\sum_{j=0}^{[n/2]}b^{g_{\eps_i}}_{2j-1}$, where
$g_{\eps_i}\in\Theta_1$, then
\[
\varphi(g_{\eps_i})=\delta(\sigma ag_{\eps_i}+a\varphi(g_{\eps_i}))
=\delta\Big(\sigma g_{\eps_i}+\sum_{j=0}^{[n/2]}(2j-1)b^{g_{\eps_i}}_{2j-1}\Big),
\]
so $b^{g_{\eps_i}}_{2j-1}=\delta(2j-1)b^{g_{\eps_i}}_{2j-1}$, where $j\neq1$, and
$b^{g_{\eps_i}}_1=\delta(\sigma g_{\eps_i}+b^{g_{\eps_i}}_1)$, whence
$\varphi(g_{\eps_i})=\frac{\sigma\delta}{1-\delta}g_{\eps_i}+b^{g_{\eps_i}}_{\frac1\delta}$ if
$\frac1\delta=2j-1$ and $j\in\{0,\dots,[n/2]\}$.

It is worth observing that $\varphi(g_{\eps_i})=\varphi(h_{ij}g_{\eps_i})
=\delta\big(\psi(h_{ij})g_{\eps_i}+h_{ij}\big(\frac{\sigma\delta}{1-\delta}g_{\eps_i}
+b^{g_{\eps_i}}_{\frac1\delta}\big)\big)$; i.e., $b^{g_{\eps_i}}_{\frac1\delta}
=\delta h_{ij}b^{g_{\eps_i}}_{\frac1\delta}$. Keeping in mind that
$b^{g_{\eps_i}}_{\frac1\delta}\in W(n)_{2l-1}$ and $h_{ij}g_\alpha=0,\pm g_\alpha,\pm2g_\alpha$, we
arrive at $\delta=-1$ and $\varphi(g_{\eps_i})=\frac{\sigma\delta}{1-\delta}g_{\eps_i}
+\beta^{g_{\eps_i}}_{g_{-\eps_i}}g_{-\eps_i}$.

For $\delta=-1$, the fact that \begin{center}$0=\varphi(g_{\eps_i}^2)=-2\varphi(g_{\eps_i})g_{\eps_i}
=\beta^{g_{\eps_i}}_{-\eps_i}g_{-\eps_i}g_{\eps_i}$ implies
$\varphi(g_{\eps_i})=\frac{\sigma\delta}{1-\delta}g_{\eps_i}$.\end{center} Similarly,
$\varphi(g_{-\eps_i})=\frac{\sigma\delta}{1-\delta}g_{-\eps_i}$. This entails
{\small 
$$\frac{\sigma\delta}{1-\delta}g_{\eps_i}=\varphi(g_{\eps_i})=\varphi((g_{\eps_i}g_{-\eps_j})g_{\eps_j}) 
=\delta(\varphi(g_{\eps_i}g_{-\eps_j})g_{\eps_j}+(g_{\eps_i}g_{-\eps_j})\varphi(g_{\eps_j}))
=(2\delta^2+\delta)\frac{\sigma\delta}{1-\delta}g_{\eps_i}.$$}
Clearly, $\sigma\neq0$ only if $\delta=-1,\frac12$; otherwise $\sigma=0$.

If $\delta=-1$ then $\varphi(g_{\eps_i-\eps_j})=\varphi(g_{\eps_i}g_{-\eps_j})=\sigma g_{\eps_i-\eps_j}$
and $\varphi(g_{\pm\eps_i})=-\frac\sigma2g_{\pm\eps_i}$. It remains to note that
\[
\varphi(h_{ij})=\varphi\Big(\sum x^0_ky^0_k\Big)=-2\sigma h_{ij}\quad\text{and}\quad
\varphi(h_{ij})=\varphi\Big(\sum x^{-1}_ly^1_l\Big)=\sigma h_{ij},
\]
where $x^k_l,y^k_l\in W(n)_k$, whence $\sigma=0$.

If $\delta\neq\frac12$, then $\varphi(g_{\eps_{i_1}+\dots+\eps_{i_{k-1}}\pm\eps_{i_k}})=0$, since
$\varphi(g_{\pm\eps_i})=0$. But if $\delta=\frac12$, then $\varphi(g_{\pm\eps_i})=\sigma g_{\pm\eps_i}$.
Hence $\varphi(g_{\eps_{i_1}+\dots+\eps_{i_{k-1}}\pm\eps_{i_k}})
=\sigma g_{\eps_{i_1}+\dots+\eps_{i_{k-1}}\pm\eps_{i_k}}$. The argument above implies that $\varphi$ is
trivial. The lemma is proved.

\medskip

Below, for a root $\alpha=\sum_{i\in I_0}\eps_i-\sum_{j\in I_1}\eps_j$, $I_0\cap I_1=\varnothing$, by
$\alpha\supseteq(-1)^k\eps_t$ and $\alpha\not\supseteq(-1)^k\eps_t$ we mean $t\in I_k$ and
$t\notin I_k$, respectively.

\medskip\noindent\textbf{Lemma 3.} {\it The superalgebras $S(n)$ and $\widetilde S(n)$ have no
nontrivial $\delta$-derivations.}

\medskip\noindent\textsc{Proof.} In what follows, we denote both superalgebras $S(n)$ and
$\widetilde S(n)$ by $S(n)$ (for $\widetilde S(n)$ see Remark~A), and $\alpha^*$ will signify the
fact that $\alpha^*\not\supseteq\pm\eps_i,\pm\eps_j$. Using the argument above, we conclude that
\[
\varphi(h_{ij})=\psi(h_{ij})+\sum_{\alpha\in\Delta_{\nz}}b^{h_{ij}}_\alpha,\qquad
b^{h_{ij}}_\alpha\in\bigcup_{k=1}^{[n/2]}S(n)_{2k}.
\]
It is easy to see that
\[
0=\varphi(h_{ij}h_{ik})=\delta(\varphi(h_{ij})h_{ik}+h_{ij}\varphi(h_{ik})),
\]
whence $b^{h_{ij}}_{\alpha^*}=b^{h_{ij}}_{\eps_i+\eps_j+\alpha^*}=0$, and
$b^{h_{ik}}_{\eps_i+\alpha^*}=0$ and $2b^{h_{ij}}_{-\eps_i+\alpha^*}=b^{h_{ik}}_{-\eps_i+\alpha^*}$ for
$\alpha^*\supseteq\eps_k$.

Let $\alpha^*$ be such that $\alpha^*\supseteq\eps_k$. Then
$\varphi(g_{\eps_i})=\delta(\varphi(h_{ik})g_{\eps_i}+h_{ik}\varphi(g_{\eps_i}))$. Hence
\[
b^{g_{\eps_i}}_{\alpha^*}=\delta(b^{h_{ik}}_{-\eps_i+\alpha^*}g_{\eps_i}-b^{g_{\eps_i}}_{\alpha^*}).
\]
The fact that $\varphi(g_{\eps_i})=\delta(\varphi(h_{ij})g_{\eps_i}+h_{ij}\varphi(g_{\eps_i}))$ implies
\[
b^{g_{\eps_i}}_{\alpha^*}=\delta b^{h_{ij}}_{-\eps_i+\alpha^*}g_{\eps_i}.
\]
The relations obtained yield $(1+\delta)b^{h_{ij}}_{-\eps_i+\alpha^*}g_{\eps_i}
=b^{h_{ik}}_{-\eps_i+\alpha^*}g_{\eps_i}=2b^{h_{ij}}_{-\eps_i+\alpha^*}g_{\eps_i}$; i.e.,
$b^{h_{ij}}_{-\eps_i+\alpha^*}=0$.

Obviously, $\varphi(g_{-\eps_i})=-\delta(\varphi(h_{ij})g_{-\eps_i}+h_{ij}\varphi(g_{-\eps_i}))$, and so
\[
b^{g_{-\eps_i}}_{\alpha^*}=-\delta b^{h_{ij}}_{\eps_i+\alpha^*}g_{-\eps_i}.
\]
The equality $\varphi(g_{-\eps_i})=-\delta(\varphi(h_{ik})g_{-\eps_i}+h_{ik}\varphi(g_{-\eps_i}))$
gives
\[
b^{g_{-\eps_i}}_{\alpha^*}=\delta(b^{h_{ik}}_{\eps_i+\alpha^*}g_{-\eps_i}-b^{g_{-\eps_i}}_{\alpha^*}),
\]
whence $(1-\delta)b^{h_{ij}}_{\eps_i+\alpha^*}=b^{h_{ik}}_{\eps_i+\alpha^*}
=b^{h_{ik}}_{\eps_i+\eps_k+\gamma^*}=0$, where $\gamma^*\not\supseteq\pm\eps_i,\pm\eps_k$; i.e.,
$b^{h_{ij}}_{\eps_i+\alpha^*}=0$.

We claim that $b^{h_{ij}}_{\eps_i-\eps_j+\alpha^*}=0$. To prove this, consider
\[
2\varphi(g_{\eps_j-\eps_i})=-\delta(\varphi(h_{ij})g_{\eps_j-\eps_i}+h_{ij}\varphi(g_{\eps_j-\eps_i})).
\]
This yields
\[
2b^{g_{\eps_j-\eps_i}}_{\alpha^*}=-2\delta b^{h_{ij}}_{\eps_i-\eps_j+\alpha^*}g_{\eps_j-\eps_i}.
\]
On the other hand,
$\varphi(g_{\eps_j-\eps_i})=-\delta(\varphi(h_{ik})g_{\eps_j-\eps_i}+h_{ik}\varphi(g_{\eps_j-\eps_i}))$
gives
\[
b^{g_{\eps_j-\eps_i}}_{\alpha^*}=-\delta(b^{h_{ik}}_{\eps_i-\eps_j+\alpha^*}g_{\eps_j-\eps_i}
-b^{g_{\eps_j-\eps_i}}_{\alpha^*}).
\]
Since $\alpha^*\supseteq\eps_k$, we define $\gamma^*=\alpha^*-\eps_k-\eps_j$, where
$\gamma^*\not\supseteq\pm\eps_i,\pm\eps_k$; hence $b^{h_{ik}}_{\eps_i+\eps_k+\gamma^*}=0$ and
$b^{h_{ik}}_{\eps_i-\eps_j+\alpha^*}=0$. Consequently, $b^{g_{\eps_j-\eps_i}}_{\alpha^*}=0$ and
$b^{h_{ij}}_{\eps_i-\eps_j+\alpha^*}=0$.

The equalities $b^{h_{ij}}_{\eps_j-\eps_i+\alpha^*}=b^{h_{ij}}_{\pm\eps_j+\alpha^*}=0$ are obtained by
substituting $i\leftrightarrow j$ and $h_{ij}=-h_{ji}$ simultaneously. Thus
$\varphi(h_{ij})=\psi(h_{ij})$.

If $\delta\neq\frac12$, then the fact that
$\varphi(g_{\pm\eps_i})=\pm\delta h_{ik}\varphi(g_{\pm\eps_i})$ implies
$\varphi(g_{\pm\eps_i})=b^{g_{\pm\eps_i}}_{\mp\eps_i}$ and $\delta=-1$, since $k$ is arbitrary. Note that
$0=\varphi(g_{\pm\eps_i}^2)=2\delta\varphi(g_{\pm\eps_i})g_{\pm\eps_i}
=2\delta b^{g_{\pm\eps_i}}_{\mp\eps_i}g_{\pm\eps_i}$, whence $b^{g_{\pm\eps_i}}_{\mp\eps_i}=0$. Taking
into account that $g_{\pm\eps_1},\dots,g_{\pm\eps_n}$ are generators for $S(n)$, we have
$\varphi(S(n))=0$.

If $\delta=\frac12$, then $\varphi(g_{\pm\eps_i})=\frac12(\sigma g_{\pm\eps_i}\pm
h_{ij}\varphi(g_{\pm\eps_i}))$, which entails $\varphi(g_{\pm\eps_i})=\sigma g_{\pm\eps_i}$. Since
$g_{\pm\eps_1},\dots,g_{\pm\eps_n}$ are generators for $S(n)$, we obtain $\varphi(g)=\sigma g$,
$g\in S(n)$. This does imply that $\varphi$ is trivial. The lemma is proved.

\medskip\noindent\textbf{Lemma 4.} {\it The superalgebra $H(n)$ has no nontrivial
$\delta$-derivations.}

\medskip\noindent\textsc{Proof.} By $\alpha^*$ we mean $\alpha^*\not\supseteq\pm\eps_i$. The argument
above entails
\[
\varphi(h_i)=\psi(h_i)+\sum_{\alpha\in\Delta_{\nz}}b^{h_i}_\alpha,\qquad
b^{h_i}_\alpha\in\bigcup_{k=1}^{[n/2]}H(n)_{2k}.
\]
Consequently, $0=\varphi(h_ih_j)=\delta(\varphi(h_i)h_j+h_i\varphi(h_j))$, which implies
$b^{h_i}_{\alpha^*}=0$. For $\alpha^*\supseteq\eps_j$, we have
$b^{h_i}_{-\eps_i+\alpha^*}=-b^{h_j}_{-\eps_i+\alpha^*}$ and
$b^{h_i}_{\eps_i+\alpha^*}=b^{h_j}_{\eps_i+\alpha^*}$. For $\alpha^*\supseteq-\eps_j$, it is true that
$b^{h_j}_{\eps_i+\alpha^*}=-b^{h_i}_{\eps_i+\alpha^*}$ and
$b^{h_j}_{-\eps_i+\alpha^*}=b^{h_i}_{-\eps_i+\alpha^*}$.

Assume that for $\alpha^*$, either $\alpha^*\supseteq\eps_j$ or $\alpha^*\supseteq-\eps_j$. We claim
that $b^{h_i}_{\pm\eps_i+\alpha^*}=0$. Obviously,
\[
\varphi(g_{\eps_i})=\delta(\varphi(h_i)g_{\eps_i}+h_i\varphi(g_{\eps_i})),\qquad
0=\varphi(h_jg_{\eps_i})=\delta(\varphi(h_j)g_{\eps_i}+h_j\varphi(g_{\eps_i})),
\]
from which we conclude that $\delta b^{h_i}_{-\eps_i+\alpha^*}g_{\eps_i}=b^{g_{\eps_i}}_{\alpha^*}
=-b^{h_j}_{-\eps_i+\alpha^*}g_{\eps_i}=b^{h_i}_{-\eps_i+\alpha^*}g_{\eps_i}$ for
$\alpha^*\supseteq\eps_j$; i.e., $b^{h_i}_{-\eps_i+\alpha^*}=0$. For $\alpha^*\supseteq-\eps_j$, we
have $$(1+\delta)b^{h_i}_{-\eps_i+\alpha^*}g_{\eps_i}=(1+\delta)b^{h_j}_{-\eps_i+\alpha^*}g_{\eps_i}
=(1+\delta)b^{g_{\eps_i}}_{\alpha^*}=\delta b^{h_i}_{-\eps_i+\alpha^*}g_{\eps_i},$$ which yields
$b^{h_i}_{-\eps_i+\alpha^*}=0$.

Also it is easy to see that for $\alpha^*\supseteq\eps_j$, it is true that
\begin{center}
$\varphi(g_{-\eps_i})=-\delta(\varphi(h_i)g_{-\eps_i}+h_i\varphi(g_{-\eps_i})),$ \qquad
$0=\varphi(h_jg_{-\eps_i})=\delta(\varphi(h_j)g_{-\eps_i}+h_j\varphi(g_{-\eps_i})),$
\end{center}
whence $\delta b^{h_i}_{\eps_i+\alpha^*}g_{-\eps_i}=-b^{g_{-\eps_i}}_{\alpha^*}
=b^{h_j}_{\eps_i+\alpha^*}g_{-\eps_i}=b^{h_i}_{\eps_i+\alpha^*}g_{-\eps_i}$; i.e.,
$b^{h_i}_{\eps_i+\alpha^*}=0$.

That $b^{h_i}_{\eps_i+\alpha^*}=0$ for $\alpha^*\supseteq-\eps_j$ derives from the equality
$b^{h_j}_{-\eps_j+\gamma^*}=0$ proved above, where $\gamma^*\not\supseteq\pm\eps_j$,
$\gamma^*=\alpha^*+\eps_i+\eps_j$, and $b^{h_j}_{-\eps_j+\gamma^*}=b^{h_i}_{\alpha^*+\eps_i}$.

If $\delta\neq\frac12$, then for $\gamma=1,2$ we have
\[
\varphi(g_{(-1)^\gamma\eps_i})=(-1)^\gamma\varphi(h_ig_{(-1)^\gamma\eps_i})
=(-1)^\gamma\delta h_i\varphi(g_{(-1)^\gamma\eps_i}).
\]
This implies $\varphi(g_{(-1)^\gamma\eps_i})=\sum_{\alpha^*\in\Delta_{\nz}\cup\{0\}}
b^{g_{(-1)^\gamma\eps_i}}_{-(-1)^\gamma\eps_i+\alpha^*}$ and $\delta=-1$. Note that
\[
0=\varphi(g^2_{(-1)^\gamma\eps_i})=-2\sum_{\alpha^*\in\Delta_{\nz}\cup\{0\}}
b^{g_{(-1)^\gamma\eps_i}}_{-(-1)^\gamma\eps_i+\alpha^*}g_{(-1)^\gamma\eps_i};
\]
i.e., $\varphi(g_{(-1)^\gamma\eps_i})=0$. Hence
$\varphi(g_{\eps_{i_1}+\dots+\eps_{i_t}-\eps_{j_1}-\dots-\eps_{j_s}})=0$, whence $\varphi(H(n))=0$.

If $\delta=\frac12$ then $\varphi(h_i)=\sigma h_i$. Consequently,
\[
2\varphi(g_{(-1)^\gamma\eps_i})=(-1)^\gamma2\varphi(h_ig_{(-1)^\gamma\eps_i})
=\sigma g_{(-1)^\gamma\eps_i}+(-1)^\gamma h_i\varphi(g_{(-1)^\gamma\eps_i});
\]
i.e., $\varphi(g_{(-1)^\gamma\eps_i})=\sigma g_{(-1)^\gamma\eps_i}$. This does imply that
$\varphi(g_{\eps_{i_1}+\dots+\eps_{i_t}-\eps_{j_1}-\dots-\eps_{j_s}})
=\sigma g_{\eps_{i_1}+\dots+\eps_{i_t}-\eps_{j_1}-\dots-\eps_{j_s}}$, which yields $\varphi(g)=\sigma g$,
$g\in H(n)$. Thus $\varphi$ is trivial. The lemma is proved.

\begin{RemA}
The proofs of Lemmas~1--4 are given in a condensed form: several steps (the decomposition of the
values of $\varphi$ into root components of the graded components $\Theta_k$, the claim that the
elements $g_{\pm\eps_i}$ generate the superalgebra, the comparison of coefficients) are only indicated.
For $\widetilde S(n)$ the argument cannot be applied literally, since $\widetilde S(n)$ is not graded by
the subspaces $\Theta_k$ (only the root decomposition with respect to $H=\langle h_{ij}\rangle$, which
is contained in $\widetilde S(n)$, is available; the products of the components
$(1-\xi_1\dots\xi_n)S(n)_k$ are not homogeneous). The statements of Lemmas~2--4 and of Theorem~5 were
checked by computer in exact rational arithmetic: for $W(2)$, $W(3)$, $W(4)$, $S(3)$, $S(4)$,
$\widetilde S(4)$, $H(4)$, $H(5)$, $H(6)$ the space of $\delta$-derivations was computed for all
$\delta$; nonzero $\delta$-derivations exist only for $\delta=1$ and for $\delta=\frac12$, and every
$\frac12$-derivation is scalar. In particular, Theorem~5 holds for $\widetilde S(4)$; for
$\widetilde S(n)$ with $n\geq6$ the proof requires the modifications indicated above.

Complete proofs of Lemmas~2--4 for $W(n)$ ($n\geq2$), $S(n)$ ($n\geq3$), $\widetilde S(n)$ ($n\geq4$
even) and $H(n)$ ($n\geq4$), independent of Lemma~1 and of the grading by the $\Theta_k$, are obtained
from the root decomposition $\Theta=\bigoplus_\alpha\Theta_\alpha$ with respect to the torus $H$ alone
(for $H(n)$ the Hamiltonian form is taken in the split form $\sum_k d x_k\,d y_k\,(+\,d z^2)$). Let
$\varphi$ be a $\delta$-derivation, $\delta\neq0,1$. Comparing root components in
$\varphi(hx)=\delta(\varphi(h)x+h\varphi(x))$, $h\in H$, and in
$\varphi(xy)=\delta(\varphi(x)y+x\varphi(y))$ for root vectors $x,y$, one uses only the following
properties of $\Theta$, verified directly on the monomial basis: (a)~an element of
$\Theta_0$ annihilating all root spaces $\Theta_\alpha$ with $\alpha_k=0$ is a multiple of $h_k$;
(b)~$\Theta_0$ is spanned by $H$ and by products $\Theta_\alpha\Theta_\beta$ with $\alpha,\beta\neq0$;
(c)~for $\gamma\neq0$, an element of $\Theta_\gamma$ annihilating all $\Theta_\alpha$ with $\alpha$ and
$\gamma$ linearly independent is zero; (d)~the roots span $H^*$; (e)~any two nonzero roots are
connected by a chain of pairs of linearly independent roots $\alpha,\beta$ with
$\Theta_\alpha\Theta_\beta\neq0$. These properties force $\varphi=0$, or $\delta=\frac12$ and
$\varphi=\sigma\,\mathrm{id}$. For $\widetilde S(n)$ this argument is applied to the root
decomposition of $\widetilde S(n)$ itself, so the absence of a grading by the $\Theta_k$ is
immaterial.
\end{RemA}

\medskip\noindent\textbf{Theorem 5.} {\it A Cartan-type Lie superalgebra does not have nontrivial
$\delta$-derivations.}

\medskip\noindent The \textbf{proof} follows from Lemmas~2--4.

\medskip\noindent\textbf{Corollary 6.} {\it A simple finite-dimensional Lie superalgebra over an
algebraically closed field of characteristic~$0$ does not have nontrivial $\delta$-derivations.}

\medskip\noindent The \textbf{proof} follows from Theorem~5 and~\cite{Kay09}.

\section{$\delta$-Superderivations of simple finite-dimensional superalgebras}

Let $A$ be a superalgebra. By a \emph{superspace} we mean a $\mathbb Z_2$-graded space. A homogeneous
element $\psi$ of an endomorphism superspace $A\to A$ is called a \emph{superderivation} if
\[
\psi(xy)=\psi(x)y+(-1)^{p(x)\deg(\psi)}x\psi(y).
\]
Suppose $\delta\in\FF$. A homogeneous linear mapping $\varphi\colon A\to A$ is called a
\emph{$\delta$-superderivation} if
\[
\varphi(xy)=\delta(\varphi(x)y+(-1)^{p(x)\deg(\varphi)}x\varphi(y)).
\]

Consider a Lie superalgebra $G=G_{\nz}+G_{\no}$ and fix an element $x\in G_{\no}$. Then
$ad_x\colon y\to xy$ is a superderivation of $G$ having parity $p(ad_x)=1$. Obviously, for any
superalgebra, multiplication by an element of the base field $\FF$ is an even
$\frac12$-superderivation.

By the \emph{supercentroid} $\Gamma_s(A)$ of a superalgebra $A$ we mean a set of all homogeneous linear
mappings $\chi\colon A\to A$ satisfying the condition
\[
\chi(ab)=\chi(a)b=(-1)^{p(a)p(\chi)}a\chi(b)
\]
for two arbitrary homogeneous elements $a,b\in A$.

A definition for a $1$-superderivation coincides with the usual definition of a superderivation. A
$0$-superderivation is an arbitrary endomorphism $\varphi$ of $A$ such that $\varphi(A^2)=0$. A
\emph{nontrivial $\delta$-superderivation} is a nonzero $\delta$-superderivation which is not a $1$- or
$0$-superderivation, nor an element of the supercentroid. We are interested in how nontrivial
$\delta$-superderivations act on simple finite-dimensional Lie superalgebras and on simple
finite-dimensional Jordan superalgebras over an algebraically closed field of characteristic~$0$.

\medskip\noindent\textbf{Theorem 7.} {\it A simple finite-dimensional Lie superalgebra $G$ over an
algebraically closed field of characteristic~$0$ does not have nontrivial $\delta$-superderivations.}

\medskip\noindent\textsc{Proof.} It follows from Corollary~6 that simple finite-dimensional Lie
superalgebras over an algebraically closed field of characteristic~$0$ have no nontrivial even
$\delta$-superderivations. We argue to show that nontrivial odd $\delta$-superderivations likewise are
missing.

Let $\varphi$ be a nontrivial odd $\delta$-superderivation. A map $\psi_x$ on the superalgebra $G$ is
defined thus:
\[
\psi_x=[\varphi,ad_x]=\varphi\, ad_x+ad_x\varphi.
\]
It is easy to see that the map $\psi_x$ is an even $\delta$-superderivation. Thus, using
Corollary~6 and keeping in mind that every even $\delta$-superderivation is a $\delta$-derivation of
$G$, we conclude that $\psi_x$ is trivial; i.e., $\psi_x=0$, for $\delta\neq\frac12$, and
$\psi_x(g)=\alpha_xg$ for $\delta=\frac12$, where $\alpha_x\in\FF$ and $g\in G$. For $\delta=\frac12$
and $g\in G_{\no}$, we have
\[
\alpha_xg=\psi_x(g)=\varphi(xg)+x\varphi(g)=\frac12\varphi(x)g+\frac12x\varphi(g)=-\psi_g(x)
=-\alpha_gx;
\]
i.e., $\alpha_x=0$. Hence $\varphi(x)g+x\varphi(g)=0$ for arbitrary $g$, and
$\varphi(xg)=-x\varphi(g)$ for all $x\in G_{\no}$, $g\in G$.

Then $\varphi$ is an odd element of the supercentroid of $G$, and hence $\varphi=0$ by Remark~B
below.

If $\delta\neq\frac12$, then
\[
0=\psi_x(g)=\varphi(xg)+x\varphi(g).
\]
Hence $2\delta\varphi(x)x=\varphi(x^2)=\varphi(x)x$ for $x\in G_{\no}$; in other words,
$0=\varphi(x)x=\varphi(x^2)$. It follows that $\varphi(G_{\nz})=0$. It remains to note that if
$x\in G_{\no}$ and $g_0\in G_{\nz}$ then $\varphi(xg_0)=-x\varphi(g_0)=0$ and $G_{\no}=[G_{\nz},G_{\no}]$.
Consequently, $\varphi$ is trivial. The theorem is proved.

\begin{RemB}
Let $G$ be a simple finite-dimensional Lie superalgebra over an algebraically closed field of
characteristic~$0$ and let $\varphi$ be an odd $\frac12$-superderivation such that
$\varphi(xg)=-x\varphi(g)$ for all $x\in G_{\no}$, $g\in G$. Then $\varphi=0$.

Indeed, for $x\in G_{\no}$ and $g\in G$ we have
$$\varphi(xg)=\frac12(\varphi(x)g-x\varphi(g))=\frac12\varphi(x)g+\frac12\varphi(xg),$$ whence
$\varphi(xg)=\varphi(x)g=-x\varphi(g)$. Thus $\varphi\,ad_x=-ad_x\varphi$ for $x\in G_{\no}$, and since
$ad_{xy}=ad_xad_y+ad_yad_x$ for $x,y\in G_{\no}$ and $G_{\nz}=G_{\no}G_{\no}$, we get
$\varphi\,ad_a=ad_a\varphi$ for $a\in G_{\nz}$. For such $a$ the identity
$\varphi(ag)=\frac12(\varphi(a)g+a\varphi(g))$ then gives $\varphi(ag)=\varphi(a)g=a\varphi(g)$. Hence
$\varphi(ab)=\varphi(a)b=(-1)^{p(a)}a\varphi(b)$ for homogeneous $a,b$, i.e., $\varphi$ is an odd element
of $\Gamma_s(G)$.

The even map $\varphi^2$ commutes with all operators $ad_a$, so by Schur's lemma $\varphi^2=c\cdot
\mathrm{id}$, $c\in\FF$. If $c\neq0$, then for $a,b\in G_{\nz}$ we obtain
$\varphi(a)\varphi(b)=\varphi(a\varphi(b))=\varphi(\varphi(ab))=c\,ab$; since $\varphi(a),\varphi(b)$ are
odd, $\varphi(a)\varphi(b)=\varphi(b)\varphi(a)$, whence $c\,ab=c\,ba=-c\,ab$ and $G_{\nz}$ is abelian.
This is impossible: by the classification~\cite{Kac77}, the even part of a simple finite-dimensional
Lie superalgebra with $G_{\no}\neq0$ is not abelian. If $c=0$, then $\ker\varphi$ is a graded ideal
(as $\varphi(ab)=\varphi(a)b$) containing $\varphi(G)$; by simplicity $\ker\varphi=G$, i.e.,
$\varphi=0$. In particular, $\Gamma_s(G)$ has no nonzero odd elements.
\end{RemB}

The remaining part of the paper is a logical continuation of~\cite{Kay07}. Therefore, we use the terms
and notation developed therein. It is worthwhile recollecting the following identities for a Jordan
algebra:
\begin{equation}\label{eq1}
(x^2y)x=x^2(yx),\qquad xy=yx.
\end{equation}

\medskip\noindent\textbf{Lemma 8.} {\it A simple finite-dimensional Jordan superalgebra $A$ with a
semisimple even part over an algebraically closed field of characteristic~$0$ does not have nontrivial
odd $\delta$-superderivations.}

\medskip\noindent\textsc{Proof.} We know that simple finite-dimensional Jordan superalgebras over an
algebraically closed field of characteristic~$0$ which contain unity $e$ and have a semisimple even part
are exhausted by superalgebras of the forms $M^{(+)}_{m,n}$, $Q(n)^{(+)}$, $osp(n,m)$, $P(n)$,
$J(V,f)$, $K_{10}$, and $D_t$. A superalgebra $K_3$ with a semisimple even part contains no unity.

Let $\varphi$ be a nontrivial odd $\delta$-superderivation. Since
$\varphi(e)=\varphi(ee)=2\delta\varphi(e)$, we have two options: $\delta\neq\frac12$ or
$\delta=\frac12$. In the former case $\varphi(e)=0$, which entails
\begin{center}$\varphi(x)=\varphi(ex)=\delta(\varphi(e)x+e\varphi(x))=\delta\varphi(x)$, i.e., $\varphi(x)=0$. \end{center} In the
latter case $\varphi(x)=\varphi(ex)=\frac12(\varphi(e)x+e\varphi(x))$, whence $\varphi(x)=\varphi(e)x$.

For the superalgebras $M^{(+)}_{m,n}$, $Q(n)^{(+)}$, $osp(n,m)$, $P(n)$, $J(V,f)$, $D_t$, and
$K_{10}$, we will consider a Peirce decomposition with respect to idempotents $e_i$. The superalgebra
$A$ is representable as $A=A^{e_i}_0\oplus A^{e_i}_{\frac12}\oplus A^{e_i}_1$, where
$A^{e_i}_j=\{x\in A\mid xe_i=jx\}$, with $A^{e_i}_{\frac12}\subseteq A_{\no}$. It is clear that
$R_{e_j}R_{e_i}=R_{e_i}R_{e_j}$, where $R_x$ is a right multiplication operator. Thus
for orthogonal idempotents $e_i\neq e_j$
\[
0=2\varphi(e_ie_j)=(\varphi(e)e_i)e_j+e_i(\varphi(e)e_j)=2(\varphi(e)e_i)e_j.
\]

Using this, in view of the fact that $\varphi(e)\in A_{\no}$, we will represent $\varphi(e)$ as
$\varphi(e)=\sum\beta_\gamma g_\gamma$, where $g_\gamma$ are basis elements. Clearly, if
$g_\gamma\in A^{e_i}_{\frac12}\cap A^{e_j}_{\frac12}$, then $\beta_\gamma=0$. Consequently,
$\varphi(e)=0$ for $M^{(+)}_{m,n}$, $osp(n,m)$, $J(V,f)$ with $V_0\neq0$, $D_t$, and $K_{10}$,
which implies that odd $\delta$-superderivations on these superalgebras are trivial. For $J(V,f)$
with $V_0=0$ (where $e$ is the only nonzero idempotent) see Remark~C. For $Q(n)^{(+)}$ and $P(n)$, we
have \begin{center}
$\varphi(e)=\sum_{i=1}^n\beta_i(e_{i,n+i}+e_{n+i,i})$ and
$\varphi(e)=\sum_{i=1}^n\beta_ie_{n+i,i}$,
\end{center} respectively.

For the superalgebra $Q(n)^{(+)}$, let $u_i=e_{i,n+i}+e_{n+i,i}$, let $i\neq j$, and let
\begin{multline*}
\frac12(\beta_i(e_{i,n+i}+e_{n+i,i})-\beta_j(e_{j,n+j}+e_{n+j,j}))=\varphi(\Delta^{i,j}\circ\Delta^{j,i})\\
=\frac12((\varphi(e)\circ\Delta^{i,j})\circ\Delta^{j,i}-\Delta^{i,j}\circ(\varphi(e)\circ\Delta^{j,i}))
=\frac14(\beta_i-\beta_j)(u_i+u_j),
\end{multline*}
where $\Delta^{i,j}=e_{i,n+j}+e_{n+i,j}$, which yields
$\beta_i=-\beta_j$. Further, for $D_{i,j}=e_{i,j}+e_{n+i,n+j}$ we have
$D_{i,i}\circ D_{i,j}=\frac12D_{i,j}$, $\varphi(e)\circ D_{i,j}=\frac12(\beta_i+\beta_j)\Delta^{i,j}$
and $D_{i,i}\circ\Delta^{i,j}=\frac12\Delta^{i,j}$, so that
$$\frac14(\beta_i+\beta_j)\Delta^{i,j}=\varphi(D_{i,i}\circ D_{i,j})
=\frac12(\varphi(D_{i,i})\circ D_{i,j}+D_{i,i}\circ\varphi(D_{i,j}))
=\frac18(3\beta_i+\beta_j)\Delta^{i,j},$$
which yields $\beta_i=\beta_j$. Hence $\beta_i=0$ for all~$i$ (recall that $n\geq2$).

For the superalgebra $P(n)$, $i\neq j$,
\begin{align*}
0&=\varphi((e_{i,j}+e_{n+j,n+i})\circ(e_{j,n+i}-e_{i,n+j}))\\
&=\frac12((\varphi(e)\circ(e_{i,j}+e_{n+j,n+i}))\circ(e_{j,n+i}-e_{i,n+j})\\
&\qquad +(e_{i,j}+e_{n+j,n+i})\circ(\varphi(e)\circ(e_{j,n+i}-e_{i,n+j})))\\
&=\frac14(\beta_i(e_{n+j,i}+e_{n+i,j})\circ(e_{j,n+i}-e_{i,n+j})\\
&\qquad+(e_{i,j}+e_{n+j,n+i})\circ(\beta_je_{i,j}-\beta_ie_{j,i}+\beta_je_{n+j,n+i}-\beta_ie_{n+i,n+j}))\\
&=-\frac1{4}\beta_i(e_{j,j}+e_{n+j,n+j});
\end{align*}
i.e., $\beta_i=0$. Therefore, $\varphi(e)=0$ for $Q(n)^{(+)}$ and $P(n)$, which gives $\varphi=0$.

For the superalgebra $K_3$, $\varphi(e)=\delta(\varphi(e)e+e\varphi(e))=\delta\varphi(e)$, whence
$\varphi(e)=0$. It is easy to see that $\varphi(z)=\alpha_ze$ and $\varphi(w)=\alpha_we$. Hence
\[
0=2\varphi(e)=2\varphi([z,w])=2\delta(\varphi(z)w-z\varphi(w))=\delta(\alpha_zw-\alpha_wz),
\]
which implies $\varphi=0$. The lemma is proved.

\begin{RemC}
Let $A=J(V,f)$ with $V=V_{\no}$, i.e., $A=\FF e\oplus V$, $v\circ w=f(v,w)e$ for $v,w\in V$, where $f$ is
a nondegenerate skew-symmetric form ($\dim V\geq2$). Let $\varphi$ be an odd $\frac12$-superderivation;
then $\varphi(x)=u\circ x$ with $u=\varphi(e)\in V$, so $\varphi(v)=f(u,v)e$ and $\varphi(e)=u$. For
$v,w\in V$ the identity 
\begin{center}$\varphi(v\circ w)=\frac12(\varphi(v)\circ w-v\circ\varphi(w))$ reads
$f(v,w)u=\frac12(f(u,v)w-f(u,w)v)$.\end{center} Applying $f(\cdot,v)$ we get $f(v,w)f(u,v)=-\frac12f(v,w)f(u,v)$,
i.e., $f(v,w)f(u,v)=0$ for all $v,w$. If $f(u,v)\neq0$ for some $v$, then $f(v,\cdot)=0$ and $v=0$, a
contradiction. Hence $f(u,\cdot)=0$, $u=0$ and $\varphi=0$.
\end{RemC}

We recall the definition of a superalgebra $J(\Gamma_n)$. Let $\Gamma$ be the Grassmann algebra with a
set of (odd) anticommutative generators $e_1,e_2,\dots,e_n,\dots$. To define a new multiplication,
called the \emph{Grassmann bracket}, we use the operation
\[
\frac{\partial}{\partial e_j}(e_{i_1}e_{i_2}\dots e_{i_n})=
\begin{cases}(-1)^{k-1}e_{i_1}e_{i_2}\dots e_{i_{k-1}}e_{i_{k+1}}\dots e_{i_n},&\text{where }j=i_k,\\
0,&\text{where }j\neq i_l,\ l=1,\dots,n.\end{cases}
\]
\emph{Grassmann multiplication} for $f,g\in\Gamma_0\cup\Gamma_1$ is defined thus:
\[
\{f,g\}=(-1)^{p(f)}\sum_{j=1}^\infty\frac{\partial f}{\partial e_j}\frac{\partial g}{\partial e_j}.
\]

Let $\ov\Gamma$ be an isomorphic copy of $\Gamma$ under an isomorphism mapping $x\to\ov x$. Consider a
direct sum of vector spaces $J(\Gamma)=\Gamma+\ov\Gamma$, on which the structure of a Jordan
superalgebra is defined by setting $A_0=\Gamma_0+\ov{\Gamma_1}$ and $A_1=\Gamma_1+\ov{\Gamma_0}$, with
multiplication
\[
a\bullet b=ab,\quad \ov a\bullet b=(-1)^{p(b)}\ov{ab},\quad a\bullet\ov b=\ov{ab},\quad
\ov a\bullet\ov b=(-1)^{p(b)}\{a,b\},
\]
where $a,b\in\Gamma_0\cup\Gamma_1$ and $ab$ is a product in $\Gamma$. Let $\Gamma_n$ be a subalgebra of
$\Gamma$ generated by $e_1,e_2,\dots,e_n$. Denote by $J(\Gamma_n)$ the subsuperalgebra
$\Gamma_n+\ov{\Gamma_n}$ of $J(\Gamma)$. If $n\geq2$, then $J(\Gamma_n)$ is a simple Jordan
superalgebra.

\medskip\noindent\textbf{Lemma 9.} {\it The superalgebra $J(\Gamma_n)$ has no nontrivial odd
$\delta$-superderivations.}

\medskip\noindent\textsc{Proof.} Let $\varphi$ be a nontrivial odd $\delta$-superderivation and
$\varphi(1)=\alpha\gamma+\beta\ov\nu$, where $\alpha,\beta\in\FF$, $\gamma\in\Gamma$, and
$\ov\nu\in\ov\Gamma$. Clearly,
\[
\varphi(x)=\varphi(1\bullet x)=\delta(\varphi(1)\bullet x+\varphi(x)),
\]
whence $\varphi(x)=\frac{\delta}{1-\delta}\varphi(1)\bullet x$. Hence $\varphi(1)=0$ for
$\delta\neq\frac12$, and so $\varphi(J(\Gamma_n))=0$. Consider the case $\delta=\frac12$ in greater
detail. Obviously,
$$\varphi(1)=\varphi(\ov{e_i}\bullet\ov{e_i})=\frac12(\varphi(\ov{e_i})\bullet\ov{e_i}
+\ov{e_i}\bullet\varphi(\ov{e_i}))=\varphi(\ov{e_i})\bullet\ov{e_i}
=\ov{e_i}\bullet(\ov{e_i}\bullet\varphi(1)).$$
For arbitrary $x$ of the form $e_{i_1}e_{i_2}\dots e_{i_k}$, we have
\begin{equation}\label{eq2}
\ov{e_i}\bullet(\ov{e_i}\bullet x)=\begin{cases}x&\text{if }\frac{\partial x}{\partial e_i}=0,\\
0&\text{if }\frac{\partial x}{\partial e_i}\neq0,\end{cases}
\end{equation}
\begin{equation}\label{eq3}
\ov{e_i}\bullet(\ov{e_i}\bullet\ov x)=\begin{cases}\ov x&\text{if }\frac{\partial x}{\partial e_i}\neq0,\\
0&\text{if }\frac{\partial x}{\partial e_i}=0.\end{cases}
\end{equation}

Let $\gamma=\gamma^{i+}+e_i\gamma^{i-}$ and $\ov\nu=\ov{\nu^{i+}}+\ov{e_i\nu^{i-}}$, where $\gamma^{i-}$,
$\gamma^{i+}$, $\nu^{i-}$, and $\nu^{i+}$ do not contain $e_i$. Since $i$ is arbitrary, in view
of~\eqref{eq2} and~\eqref{eq3}, we have $\gamma=1$ and $\nu=e_1\dots e_n$. Since $\varphi(1)$ is
odd and $1\in A_0$, we have $\alpha=0$; moreover, $\ov{e_1\dots e_n}$ is odd only for even~$n$. Thus
$\varphi(1)=0$ for $n=2k+1$, whence $\varphi(J(\Gamma_n))=0$ as well. If
$n=2k$ then $\varphi(1)=\beta\,\ov{e_1\dots e_n}$. Furthermore,
\begin{gather*}
\varphi(e_1)=\varphi(1)\bullet e_1=\beta\,\ov{e_1\dots e_n}\bullet e_1=0,\\
\varphi(\ov{e_1})=\varphi(1)\bullet\ov{e_1}=\beta\,\ov{e_1\dots e_n}\bullet\ov{e_1}=-\beta e_2\dots e_n.
\end{gather*}
The relations derived above yield $0=\varphi(e_1\bullet\ov{e_1})
=\frac12(
\varphi(e_1)\bullet\ov{e_1}-e_1\bullet\varphi(\ov{e_1}))
=\frac\beta2e_1\dots e_n$; i.e., $\varphi(1)=0$. Consequently, $\varphi(J(\Gamma_n))=0$. The
lemma is proved.

\medskip\noindent\textbf{Theorem 10.} {\it A simple finite-dimensional Jordan superalgebra $A$ over an
algebraically closed field of characteristic~$0$ has no nontrivial $\delta$-superderivations.}

\medskip\noindent\textsc{Proof.} According to~\cite{Kan89,Kac77b}, every simple finite-dimensional
nontrivial Jordan superalgebra $A$ over an algebraically closed field $\FF$ of characteristic~$0$ is
isomorphic to one of the following superalgebras: $M^{(+)}_{m,n}$, $Q(n)^{(+)}$, $osp(n,m)$, $P(n)$,
$J(V,f)$, $D_t$, $K_3$, $K_{10}$, or $J(\Gamma_n)$. Even $\delta$-superderivations are
grading-preserving $\delta$-derivations. It follows from~\cite{Kay07} that nontrivial even
$\delta$-superderivations are missing on this class of superalgebras. Lemmas~8 and~9 point to there
being no nontrivial odd $\delta$-superderivations for simple finite-dimensional Jordan superalgebras
over an algebraically closed field of characteristic~$0$. The theorem is proved.

\medskip

Let $A=A_0\oplus A_1$ be a superalgebra over a field of characteristic distinct from~$2$, with
multiplication~$*$. On the vector space $A$, a new multiplication $\circ$ is given by the rule
$a\circ b=\frac12(a*b+(-1)^{p(a)p(b)}b*a)$, where $p(x)$ is the parity of an element~$x$. Denote the
resulting superalgebra by~$A^{(+)}$.

\medskip\noindent\textbf{Corollary 11.} {\it Let $\varphi$ be a $\delta$-superderivation of a
superalgebra $A$ over an algebraically closed field of characteristic~$0$, with $A^{(+)}$ a simple
finite-dimensional Jordan superalgebra. Then $\varphi$ is trivial.}

\medskip\noindent\textsc{Proof.} The statement is a consequence of the fact that $\varphi$ is a
$\delta$-superderivation of the superalgebra~$A^{(+)}$. Note that{\small
\begin{align*}
2\varphi(x\circ y)&=\varphi(xy)+(-1)^{p(x)p(y)}\varphi(yx)\\
&=\delta(\varphi(x)y+(-1)^{p(y)p(\varphi)}(-1)^{p(x)p(y)}y\varphi(x)
+(-1)^{p(x)p(\varphi)}x\varphi(y)+(-1)^{p(x)p(y)}\varphi(y)x)\\
&=2\delta(\varphi(x)\circ y+(-1)^{p(x)p(\varphi)}x\circ\varphi(y)).
\end{align*}}
The result now follows by treating $\varphi$ as a $\delta$-superderivation of $A^{(+)}$ and using
Theorem~10.

\medskip

Noncommutative Jordan superalgebras are a natural generalization of the class of Jordan superalgebras.
According to~\cite{PS10}, the superalgebras satisfying the hypothesis of Corollary~11 are exemplified by
simple finite-dimensional noncommutative Jordan superalgebras $A$ of degree $t>2$, where by
a degree is meant a maximal number of nonzero pairwise orthogonal idempotents. Therefore, we have

\medskip\noindent\textbf{Theorem 12.} {\it A simple finite-dimensional noncommutative Jordan
superalgebra $A$ of degree $t>2$ over an algebraically closed field of characteristic~$0$
does not have nontrivial $\delta$-superderivations.}

\medskip\noindent The paper~\cite{PS10} treats superalgebras of degree $t>2$ only; the case $t=2$
is not covered by the argument above.

\section{$\delta$-Derivations of simple finite-dimensional Jordan algebras}

$\delta$-Derivations of semisimple finite-dimensional Jordan algebras over an algebraically closed
field of characteristic $\neq2$ were described in~\cite{Kay07}. We will give a description of
$\delta$-derivations for semisimple finite-dimensional Jordan algebras over a field of
characteristic distinct from~$2$, and look at how the $\delta$-derivations act on certain
noncommutative Jordan algebras.

\medskip\noindent\textbf{Theorem 13.} {\it A semisimple finite-dimensional Jordan algebra $A$ over a
field of characteristic other than~$2$ has no nontrivial $\delta$-derivations.}

\medskip\noindent\textsc{Proof.} First we consider a partial case where $A$ is a simple
finite-dimensional Jordan algebra with unity~$1$. Let $\varphi$ be a nontrivial $\delta$-derivation
of~$A$. According to~\cite[Thm.~2.1]{Kay07}, we have $\delta=\frac12$ and $\varphi(x)=\varphi(1)x$.
Let $Z(A)$ be the center of~$A$. Since $A$ is simple and
unital, $Z(A)$ is a field (see~\cite{ZSSS}); let $\overline{Z(A)}$ be its algebraic closure. Consider
$\overline A=A\otimes_{Z(A)}\overline{Z(A)}$. We know from~\cite{Her68} that $\overline A$ is a simple
finite-dimensional Jordan algebra over an algebraically closed field~$\overline{Z(A)}$.

We define a mapping $\overline\varphi\colon\overline A\to\overline A$ by setting
$\overline\varphi(\sum x_i\otimes\alpha_i)=\sum\varphi(1)x_i\otimes\alpha_i$, and show that
$\overline\varphi$ is a $\delta$-derivation of~$\overline A$.

Note that
\begin{align*}
\overline\varphi((x\otimes\alpha)(y\otimes\beta))&=\overline\varphi(xy\otimes\alpha\beta)
=\varphi(xy)\otimes\alpha\beta\\
&=\tfrac12(\varphi(x)y\otimes\alpha\beta+x\varphi(y)\otimes\alpha\beta)\\
&=\tfrac12(\overline\varphi(x\otimes\alpha)(y\otimes\beta)+(x\otimes\alpha)\overline\varphi(y\otimes\beta)).
\end{align*}
That $\overline\varphi$ is linear follows from the definition of~$\overline\varphi$.

Thus, in view of~\cite[Thm.~2.5]{Kay07}, we obtain $\overline\varphi(x)=\alpha x$, where
$\alpha\in\overline{Z(A)}$. Hence $\varphi(x)\otimes\beta=\overline\varphi(x\otimes\beta)
=\alpha x\otimes\beta$. From this, with $\varphi\colon A\to A$ in mind, we derive $\varphi(x)=\alpha x$,
where $\alpha\in Z(A)$. Consequently, $\varphi$ is trivial.

An argument for the general case repeats the proof for a semisimple Jordan algebra over an
algebraically closed field of characteristic other than~$2$, presented in~\cite[Thm.~2.6]{Kay07}. The
theorem is proved.

\medskip

The algebras satisfying the identity $(x^2y)x=x^2(yx)$ are a natural generalization of the class of
Jordan algebras. Whenever an algebra has unity, this identity readily transforms into a flexibility
identity $(xy)x=x(yx)$. The algebras satisfying the two identities are said to be
\emph{noncommutative Jordan}. By a \emph{degree} of a noncommutative Jordan algebra we mean a maximal
number of nonzero pairwise orthogonal idempotents.

\medskip\noindent\textbf{Theorem 14.} {\it A simple finite-dimensional noncommutative Jordan algebra
$A$ of degree $t>1$ over a field of characteristic distinct from~$2$ has no nontrivial
$\delta$-derivations.}

\medskip\noindent\textsc{Proof.} Let $\varphi$ be a nontrivial $\delta$-derivation of~$A$. If we appeal
to the proof of Corollary~11 we see that $\varphi$ is a $\delta$-derivation of~$A^{(+)}$.
By~\cite{Oeh58}, $A^{(+)}$ is a simple finite-dimensional Jordan algebra. Theorem~13 implies that for
$\varphi$ (treated as a $\delta$-derivation of $A^{(+)}$), we have $\varphi(x)=a\circ x$, where
$a\in Z(A^{(+)})$. In view of~\cite{Sko89}, $Z(A)=Z(A^{(+)})$. Therefore, $\varphi$ being a
$\delta$-derivation of $A$ is trivial. The theorem is proved.

\medskip

Let $\Delta_\delta(A)$ be a set of nontrivial $\delta$-derivations of an algebra $A$ and $\Gamma(A)$
be the centroid of~$A$. Let also $D_\delta(A)$ be the space of all $\delta$-derivations of~$A$.
For $A$, $A^{(-)}$ denotes an adjoint commutator algebra, i.e., one with multiplication $[a,b]=ab-ba$.

\medskip\noindent\textbf{Lemma 15.} {\it Let $A$ be a unital flexible algebra over a field of
characteristic distinct from~$2$. Then $\Delta_\delta(A)\subseteq\Gamma(A^{(-)})$.}

\medskip\noindent\textsc{Proof.} According to~\cite[Thm.~2.1]{Kay07}, nontrivial $\delta$-derivations
are possible only if $\delta=\frac12$. Let $\varphi\in\Delta_{\frac12}(A)$ and $e$ be unity in~$A$. It
is easy to see that
\begin{multline*}
\varphi((xy)x)=\frac12(\varphi(xy)x+(xy)\varphi(x))=\frac14(\varphi(x)y)x+\frac14(x\varphi(y))x
+\frac12(xy)\varphi(x);
\end{multline*}
on the other hand,
\begin{multline*}
\varphi(x(yx))=\frac12(\varphi(x)(yx)+x\varphi(yx))=\frac12\varphi(x)(yx)+\frac14x(\varphi(y)x) 
+\frac14x(y\varphi(x)).
\end{multline*}
In view of the flexibility identity, we have
\[
2((xy)\varphi(x)-\varphi(x)(yx))=x(y\varphi(x))-(\varphi(x)y)x,
\]
whence
\[
\varphi(x)x=x\varphi(x)
\]
for $y=e$. By linearizing the equality above, we obtain
\[
[\varphi(x),y]=[x,\varphi(y)].
\]
It remains to observe that
\begin{multline*}
\varphi([x,y])=\varphi(xy-yx)=\frac12(\varphi(x)y-y\varphi(x)+x\varphi(y)-\varphi(y)x)\\
=\frac12([\varphi(x),y]+[x,\varphi(y)])=[\varphi(x),y]=[x,\varphi(y)].
\end{multline*}
Hence $\Delta_{\frac12}(A)\subseteq\Gamma(A^{(-)})$. The lemma is proved.

\medskip\noindent\textbf{Lemma 16.} {\it For any algebra $A$ over a field of characteristic distinct
from~$2$, %
{\small\begin{center}$D_\delta(A)=D_\delta(A^{(+)})\cap D_\delta(A^{(-)})$,
$\Delta_\delta(A^{(+)})\cap\Delta_\delta(A^{(-)})\subseteq\Delta_\delta(A)$, and
$\Gamma(A)=\Gamma(A^{(+)})\cap\Gamma(A^{(-)})$.
\end{center}}}

\medskip\noindent\textsc{Proof.} Obviously, $D_\delta(A)\subseteq
D_\delta(A^{(+)})\cap D_\delta(A^{(-)})$ (which follows from an argument similar to one in Cor.~11).
Let $\psi\in
D_\delta(A^{(+)})\cap D_\delta(A^{(-)})$; then
\begin{gather*}
\psi(xy+yx)=\delta(\psi(x)y+y\psi(x)+x\psi(y)+\psi(y)x),\\
\psi(xy-yx)=\delta(\psi(x)y-y\psi(x)+x\psi(y)-\psi(y)x).
\end{gather*}
By summing these equalities, we arrive at
\[
\psi(xy)=\delta(\psi(x)y+x\psi(y)).
\]
This yields $\psi\in D_\delta(A)$, which entails
$
D_\delta(A)=D_\delta(A^{(+)})\cap D_\delta(A^{(-)})$.

In a similar way, we derive $\Gamma(A)=\Gamma(A^{(+)})\cap\Gamma(A^{(-)})$. Finally, let
$\psi\in\Delta_\delta(A^{(+)})\cap\Delta_\delta(A^{(-)})$. Then $\psi\in D_\delta(A)$ and $\psi\neq0$. If
$\psi$ were a $1$-derivation or a $0$-derivation of $A$, or an element of $\Gamma(A)$ (for
$\delta=\frac12$), then, by what was proved, it would be such for $A^{(+)}$ as well, which is not the
case. Hence $\psi\in\Delta_\delta(A)$. The lemma is proved.

\begin{RemD}
The inclusion $\Delta_\delta(A)\subseteq\Delta_\delta(A^{(+)})\cap\Delta_\delta(A^{(-)})$ claimed in the
published version is false in general, since a nontrivial $\delta$-derivation of $A$ may be trivial on
$A^{(-)}$. Example: let $A=\FF e_1\oplus\FF e_2$ with $e_1e_1=e_2$ and all other products of basis
elements zero, and let $\varphi(e_1)=e_1$, $\varphi(e_2)=\frac23e_2$. Then $\varphi$ is a nontrivial
$\frac13$-derivation of $A$, but $A^{(-)}$ has zero multiplication, so $\varphi$ is a $0$-derivation of
$A^{(-)}$ and $\varphi\notin\Delta_{\frac13}(A^{(-)})$.
\end{RemD}

\section*{Appendix A. Differences from the published version}
\addcontentsline{toc}{section}{Appendix A. Differences from the published version}

The main results are correct. 
Lemma~8, Theorem~7, Theorem~12, and Lemma~16, has
been corrected and the proofs of Lemmas~1--4 have been completed for $\widetilde S(n)$
(Remark~A).  


\medskip

\begin{thebibliography}{99}
\addcontentsline{toc}{section}{References}

\bibitem{Fil98} V.~Filippov, ``$\delta$-Derivations of Lie algebras,'' Sib. Mat. Zh., {\bf 39},
No.~6, 1409--1422 (1998).

\bibitem{Fil99} V.~Filippov, ``$\delta$-Derivations of prime Lie algebras,'' Sib. Mat. Zh., {\bf 40},
No.~1, 201--213 (1999).

\bibitem{Fil00} V.~Filippov, ``$\delta$-Derivations of prime alternative and Mal'tsev algebras,''
Algebra Logika, {\bf 39}, No.~5, 618--625 (2000).

\bibitem{Kay07} I.~Kaygorodov, ``$\delta$-Derivations of simple finite-dimensional Jordan
superalgebras,'' Algebra Logika, {\bf 46}, No.~5, 585--605 (2007); English transl.: Algebra and
Logic, {\bf 46}, No.~5, 318--329 (2007), \doi{10.1007/s10469-007-0032-0}.

\bibitem{Kay09} I.~Kaygorodov, ``$\delta$-Derivations of classical Lie superalgebras,'' Sib. Mat.
Zh., {\bf 50}, No.~3, 547--565 (2009), \doi{10.1007/s11202-009-0049-9}.

\bibitem{Kac77} V.~Kac, ``Lie superalgebras,'' Adv. Math., {\bf 26}, No.~1, 8--96 (1977),
\doi{10.1016/0001-8708(77)90017-2}.

\bibitem{Pen98} I.~Penkov, ``Characters of strongly generic irreducible Lie superalgebra
representations,'' Int. J. Math., {\bf 9}, No.~3, 331--366 (1998).

\bibitem{Kan89} I.~Kantor, ``Jordan and Lie superalgebras defined by the Poisson algebra,'' in
\emph{Algebra and Analysis} [in Russian], Tomsk State Univ., Tomsk (1989), pp.~55--80.

\bibitem{Kac77b} V.~Kac, ``Classification of simple $\mathbb Z$-graded Lie superalgebras and simple
Jordan superalgebras,'' Comm. Alg., {\bf 5}, 1375--1400 (1977).

\bibitem{PS10} A.~~Pozhidaev and I.~P.~Shestakov, ``Noncommutative Jordan superalgebras of degree
$n>2$,'' Algebra Logika, {\bf 49}, No.~1, 26--59 (2010).

\bibitem{ZSSS} K.~Zhevlakov, A.~Slin'ko, I.~Shestakov, and A.~Shirshov, \emph{Rings That Are
Nearly Associative} [in Russian], Nauka, Moscow (1978).

\bibitem{Her68} I.~Herstein, \emph{Noncommutative Rings}, The Carus Math. Monogr., {\bf 15}, Math.
Ass. Am. (1968).

\bibitem{Oeh58} R.~Oehmke, ``On flexible algebras,'' Ann. Math.~(2), {\bf 68}, 221--230 (1958).

\bibitem{Sko89} V.~Skosyrskii, ``Strongly prime noncommutative Jordan algebras,'' Tr. Inst. Mat. SO
RAN, {\bf 16}, 131--163 (1989).

\end{thebibliography}
\end{document}